\documentclass[12pt]{amsart}

\usepackage{amscd,amsmath}

\newcommand{\la}{\lambda}

\newcommand{\norm}[1]{\Vert #1\Vert}

\newcommand{\abs}[1]{\vert #1\vert}

\newcommand{\ud}{\underline}

\newcommand{\Ll}{\mathcal{L}}

\newcommand{\f}{\varphi}
\newcommand{\vv}{\overrightarrow}

\newcommand{\CC}{\mathbb{C}}   
   
\newcommand{\RR}{\mathbb{R}}   
\newcommand{\ZZ}{\mathbb{Z}}

\numberwithin{equation}{section}

\newtheorem{te}{Theorem}[section]

\theoremstyle{definition}
\newtheorem{de}{Definition}[section]    

\theoremstyle{remark}

\newtheorem{ex}{Example}[section]

\begin{document}

\title[Sasakian reduction]{ Reduction of Sasakian manifolds}

\author{Gueo Grantcharov}
\address{Department of Mathematics\\ University of California at Riverside\\
Riverside CA 92521, U.S.A.}
\email{geogran@math.ucr.edu}
\author{Liviu Ornea}
\address{University of Bucharest\\ Faculty of Mathematics\\ 14 Academiei str.\\ 70109
Bucharest, Romania}
\email{lornea@geometry.math.unibuc.ro}

\date{\today}

\begin{abstract}
We show that  the contact reduction  can be specialized to Sasakian 
manifolds. We link this Sasakian reduction to K\"ahler reduction 
by considering the K\"ahler cone over a Sasakian manifold. We present examples of
Sasakian manifolds obtained by $S^1$ reduction of standard Sasakian spheres.
\end{abstract} 

\maketitle

\section{Introduction}
Reduction technique was naturally extended from symplectic to contact 
structures by H. Geiges in \cite{Gei}. Even earlier, Ch. Boyer, K. Galicki 
and B. Mann defined in \cite{Bo-Ga-Ma 1} a moment map for $3$-Sasakian 
manifolds, thus extending 
the reduction procedure for nested metric contact structures. Quite surprisingly, 
a reduction scheme for Sasakian manifolds (contact manifolds endowed with a 
compatible Riemannian metric satisfying a curvature condition),  
was still missing. 

    In this note we fill the gap by defining a Sasakian moment map and constructing 
the associated reduced space. We then relate Sasakian reduction to K\"ahler 
reduction \emph{via} the K\"ahler cone over a Sasakian manifold.

	{\bf Acknowledgements} This research was initiated during the authors 
visit at the 
\emph{Abdus Salam International 
Centre for Theoretical Physics}, Trieste, during summer 1999. The authors 
thank the Institute for support and excellent environment. The second author
also aknowledges financial and technical support  
from the \emph{Erwin Schr\"odinger Institute},
Vienna, in September 1999. Both authors are grateful to Kris Galicki and Henrik
Pedersen for many 
illuminating conversations on Sasakian geometry and related themes. 

\section{Sasakian manifolds}
Let us briefly recall the notion of a Sasakian manifold. The definition we give
is not the standard one but is suited for our purpose. For more details, 
we refer to  \cite{Bl} and \cite{Bo-Ga 3} .  

\begin{de}
A {\it Sasakian} manifold is a $(2n+1)$-dimensional Riemannian manifold $(N,g)$
endowed with a unitary Killing vector field $\xi$ such that the 
curvature tensor of $g$ satisfies the equation:
\begin{equation}\label{curva}
R(X,\xi)Y=\eta(Y)X-g(X,Y)\xi
\end{equation}
 where
$\eta$ is the metric 
dual 1-form of $\xi$: $\eta(X)=g(\xi,X)$.
\end{de}
Let $\f=\nabla\xi$, whwre $\nabla$ is the Levi-Civita connection of $g$. 
The following formulae are then easily deduced:
\begin{equation}
\f \xi=0,\quad g(\f Y,\f Z)=g(Y,Z)-\eta(Y)\eta(Z).
\end{equation}
It can be seen that 
$\eta$ is a contact form on $N$, whose  Reeb field is $\xi$ (it is also called the 
characteristic vector field). Moreover, the restriction of $\f$ to the contact
distribution $\eta=0$ is a complex structure.  

     The simplest example is the standard sphere
$S^{2n+1}\subset \CC^{n+1}$, with the metric induced by the flat one of $\CC^{n+1}$. 
The characteristic  Killing vector field is $\xi_p=-i\vv{p}$, $i$ being the 
imaginary unit. Other Sasakian structures on the sphere can be obtained by
$D$-homothetic transformations (cf. \cite{Ta}). Also, the unit sphere bundle of any
space form is Sasakian. 

    More 
generally, the quantization bundle of a compact K\"ahler manifold naturally 
carries a 
Sasakian structure. The converse construction, possible when the characteristic 
field is regular, is known as the Boothby-Wang fibration. Precisely, the 
following result  (the metric part is due to Morimoto and Hatakeyama) 
is available (cf. \cite{YK} or \cite{Bo-Ga 3}):
\begin{te}
Let $(P,h)$ be a Hodge manifold. There exists a principal circle bundle 
$\pi:N\rightarrow P$ 
and a connection form $\eta $ in it, with curvature form the pull-back of 
the K\"ahler form of $P$, which is a contact form on $S$. Let $\xi$ be the  
vector field dual to $\eta$ with respect to the metric 
$g=\pi^*h+\eta\otimes\eta$. Then $(N,g,\xi)$ is Sasakian.
\end{te}

    The following equivalent definition puts Sasakian geometry in the framework 
of holonomy groups. Let $C(N)=N\times \RR_+$ be the cone over $(N,g)$. Endow 
it with the warped-product cone metric $C(g)=r^2g+dr^2$. Let $R_0=r\partial r$ 
and define on 
$C(N)$ the complex structure $J$ acting like this (with obvious 
identifications):  
$JY=\f Y-\eta(Y)R_0$, $JR_0=\xi$. We have:
\begin{te}\cite{Bo-Ga 3} 
$(N,g,\xi)$ is Sasakian if and only if the cone over $N$ 
$(C(N),C(g),J)$ is K\"ahlerian. 
\end{te}

\section{Main results}
\begin{te}
Let $(N,g,\xi)$ be a compact $2n+1$ dimensional Sa\-sa\-ki\-an manifold and $G$ a 
compact $d$-dimensional Lie group acting on $N$ by contact isometries. Suppose 
$0\in \mathfrak{g}^*$ is a regular value of the associated moment map $\mu$. Then 
the reduced space $M=N/\!/G:=\mu^{-1}(0)/G$ is a Sasakian manifold of dimension 
$2(n-d)+1$.
\end{te}
\begin{proof}
By \cite{Gei}, the contact moment map $\mu:N\rightarrow \mathfrak{g}^*$
is defined by
\begin{equation*}
<\mu(x), \underline{X}>=\eta(X)
\end{equation*}
for any $\underline{X}\in \mathfrak{g}$ and $X$ the corresponding field
on $N$. We know that the reduced space is a
contact manifold, \emph{loc. cit}. Hence we 
only need to check that (1) the Riemannian metric is projected on $M$ and 
(2) the field $\xi$ projects to a unitary Killing field on $M$ such that 
the curvature tensor of the projected metric satisfies formula \eqref{curva}. 

	To this end, we first 
describe the metric geometry of the Riemannian submanifold $\mu^{-1}(0)$.

    Let $\{\ud{X}_1,...,\ud{X}_d\}$ be a basis of $\mathfrak{g}$ and let $\{X_1,...,
X_d\}$ be the corresponding vector fields on $N$. Since $0$ is a regular value 
of $\mu$, $\{{X_i}_x\}$ is a linearly independent system in each 
$T_x\mu^{-1}(0)$. From the very definition of the 
moment map we have $\eta_p(X_i)=\mu(p)(X_i)=0$ hence $X_i\perp \xi$. 
As $G$ acts by contact isometries, we have 
\begin{equation}
\Ll_{X_i}g=0, \quad \Ll_{X_i}\eta=0\quad i=1,...,d.
\end{equation}
Note that these also imply $[X_i,\xi]=\Ll_{X_i}\xi=0$. 

Observe that $\mu^{-1}(0)$ is an 
isometrically immersed submanifold of $N$ (we denote the induced metric also 
with 
$g$) whose tangent space in 
each point is described by: $Y\in T_x\mu^{-1}(0)$ if and only if 
$d\mu_x(Y)=0$. Hence, 
by the definition of the moment map,   
the vector fields $\xi$ and 
$X_i$ are tangent to $\mu^{-1}(0)$. Moreover, for any $Y$ tangent to 
$\mu^{-1}(0)$, one has $g(\f X_i, Y)=d\eta(Y,X_i)=d\mu(Y)=0$, hence the 
vector fields $\{X_i\}$ produce a local basis (not necessarily orthogonal) of 
the normal bundle of $\mu^{-1}(0)$. The shape operators $A_i:=A_{\f X_i}$ 
 of this submanifold in $N$ are computed as follows (we let $\nabla$, $\nabla^N$ 
 be the Levi Civita covariant derivatives of $\mu^{-1}(0)$, resp. $N$):
\begin{equation}\label{2f}
\begin{split}
g(A_iY,Z)&=-g(\nabla^N_Y(\norm{X_i}^{-1}\f X_i),Z)=\\
&=-g(Y(\norm{X_i}^{-1})\f X_i,Z)-g(\norm{X_i}^{-1}\nabla^N_Y(\f X_i),Z)=\\
&=-\norm{X_i}^{-1}g(\nabla^N_Y(\f X_i),Z)=\\
&=-\norm{X_i}^{-1}g(\nabla^N_Y(\f) X_i+\f \nabla^N_YX_i,Z)=\\
&=-\norm{X_i}^{-1}g(\eta(X_i)Y-g(X_i,Y)\xi+\f \nabla^N_YX_i,Z)=\\
&=\norm{X_i}^{-1}\{g(X_i,Y)\eta(Z)-g(\f \nabla^N_YX_i,Z\}.
\end{split}
\end{equation}
In particular, for the corresponding quadratic second funadamental forms we get:
 \begin{equation}\label{unu}
h_i(Y,\xi)=\norm{X_i}^{-1}g(X_i,Y),\quad h_i(\xi,\xi)=0.
\end{equation}
Consequently, one easily obtains:  
\emph{the restriction of the vector field $\xi$ is Killing on $\mu^{-1}(0)$ 
too}. 

    Using the Gauss equation of a submanifold
\begin{equation*}
\begin{split}
R^N(X,Y,Z,W)&=
R^{\mu^{-1}(0)}(X,Y,Z,W)\\
&+g(h(X,Z),h(Y,W))-g(h(X,W),h(Y,Z))
\end{split}
\end{equation*}   
 and the formula \eqref{2f} we now compute the  needed part of
the curvature tensor of  
 $\mu^{-1}(0)$  at a fixed point 
$p\in \mu^{-1}(0)$. We take $X,Y,Z$ orthogonal to $\xi_p$ and obtain:
\begin{equation}\label{poon}
\begin{split}
g(R&^{\mu^{-1}(0)}(X,\xi)Y,Z)-g(R^N(X,\xi)Y,Z)=\\
&=-\sum_{i=1}^d\norm{X_i}^{-2}
\left\{h_i(X,Y)h_i(\xi,Z)
-h_i(X,Z)h_i(\xi,Y)\right\}\\
&=
-\sum_{i=1}^d\norm{X_i}^{-2}
\left\{g(X_i,Z)g(\nabla^N_XX_i,\f Y)-g(X_i,Y)g(\nabla^N_XX_i,\f Z)\right\}
\end{split}
\end{equation}
(Note  that $\nu_i=\norm{X_i}^{-1}\f {X_i}_p$ are chosen to be 
orthonormal in $p$;  this is always possible pointwise by 
appropriate choice of the initial $\ud{X}_i$).

    Let now $\pi:\mu^{-1}(0)\rightarrow M$ and endow $M$ with the projection $g^M$  
of the metric $g$ such that $\pi$ becomes a Riemannian 
submersion. This is possible because $G$ acts by isometries. In this setting, 
the vector fields $X_i$ span the vertical distribution of the submersion, whilst 
$\xi$ is horizontal and projectable (because $\Ll_{X_i}\xi=0$). Denote with 
$\zeta$ its projection on $M$. $\zeta$ is obviously  unitary. 
To prove that $\zeta$ is Killing on $M$, we just observe that 
$\Ll_\zeta g(Y,Z)=\Ll_\xi g(Y^h,Z^h)$, where  $Y^h$ denotes the horizontal lift of 
of $Y$. Finally, to compute the  values $R^M(X,\zeta)Y$ of the curvature tensor 
of $g^M$, we use O'Neill formula (cf. \cite{Be}, (9.28f))
\begin{equation*}
\begin{split}
g^M(&R^M(X,\zeta)Y,Z)=g(R^{\mu^{-1}(0)}(X^h,\xi)Y^h,Z^h)+g(A(X^h,\xi),A(Y^h,Z^h))\\&-
g(A(\xi,Y^h),A(X^h,Z^h))+g(A(X^h,Z^h),A(\xi,Z^h))
\end{split}
\end{equation*}
where $X,Y,Z$ are unitary, normal to $\zeta$ and the O'Neill $(1,2)$ tensor $A$ is 
defined as: 
$A(Z^h,X^h)=\text{\emph{vertical part of}}\, \nabla_{Z^h}X^h$. Using Gauss formula 
and \eqref{unu}, 
we obtain
\begin{equation*}
g(\nabla_{Z^h}\xi, X_i)=g(\f Z^h,X_i)=-g(Z^h,\f X_i)=0
\end{equation*}
hence $\nabla_{Z^h}\xi$ has no vertical part and $A(Z^h,\xi)=0$. Thus  
$$R^M(X,\zeta)Y=R^{\mu^{-1}(0)}(X^h,\xi)Y^h=R^N(X^h,\xi)Y^h$$
because of \eqref{poon} and the fact that $X^h, Y^h$ are normal to all $X_i$.
 Hence 
$$R^M(X,\zeta)Y=g(\xi,Y^h)X^h-g(X^h,y^h)\xi=g^M(\zeta,Y)X-g^M(X,Y)\zeta$$
which  proves that $(M,g^M,\zeta)$ is a Sasakian manifold.
\end{proof}
\medskip

    In the following we relate Sasakian reduction to  K\"ahler reduction by using 
the cone construction. Roughly speaking, we prove that reduction and taking 
the cone are commuting operations.

    Let $\omega=dr^2\wedge \eta+r^2d\eta$ be the K\"ahler form of the cone 
$C(N)$ over a Sasakian manifold $(N,g,\xi)$. If $\rho_t$ are the translations acting on 
$C(N)$ by $(x,r)\mapsto (x,tr)$, then the vector field $R_0=r\partial r$ is 
the one generated by $\{\rho_t\}$. Moreover, the following two relations 
are useful:
\begin{equation}
\Ll_{R_0}\omega=\omega, \quad \rho_t^*\omega=t\omega.
\end{equation}
Suppose a compact Lie group $G$ acts on $C(N)$ by holomorphic isometries, 
commuting with $\rho_t$.  This ensures 
 a corresponding action of $G$ on $N$. In fact, we can consider $G\cong G\times\{Id\}$ acting 
 as $(g, (x,r))\times (gx,r)$. 
     
     Suppose  
that a moment map $\Phi:C(N)\rightarrow \mathfrak{g}$ exists.
 
     As above, let 
$\{\ud{X}_1,...,\ud{X}_d\}$ be a basis of $\mathfrak{g}$ and let $\{X_1,...,
X_d\}$ be the corresponding vector fields on $C(N)$. We see that $X_i$ are 
independent on $r$, hence can be 
considered as vector fields on $N$. Furthermore, the commutation of $G$ with 
$\rho_t$ implies 
\begin{equation}\label{doi}
\Phi(\rho_t(p))=t\Phi(p).
\end{equation} 
Now imbed $N$ in the cone as 
$N\times \{1\}$ and let $\mu:=\Phi|_{N\times\{1\}}$. This is the moment 
map of the action of $G$ on $N$. 
        To see this, recall the definition of the symplectic moment map 
$\Phi=(\Phi_1,...,\Phi_d)$:  
$\Phi_i$  is given up to constant by $d\Phi_i(Y)=\omega(X_i,Y)$. Here we 
uniquely determine $\Phi_i$ by imposing the condition 
$\eta(X_i)=\Phi|_{N\times\{1\}}$.  This immediately implies that the Reeb field 
of $N$ is orthogonal to the vector fields $X_i$ since $g(\xi, X_i)=\eta(X_i)=0$. 
As $G$ acts by isometries on $C(N)$, we may project the cone metric to 
a metric on  $ N'/\!/G\times \RR_+$ which we denote by $g_0$. 
Then $g_0(Y,Z)=C(g)(Y^h,Z^h)$, 
where $Y^h$, $Z^h$ are the unique vector fields on $\Phi^{-1}(0)$ orthogonal to all 
of $X_i$ which project on $Y$, $Z$ (we call them horizontal).

    Let $P=\Phi^{-1}(0)/G$ be the reduced K\"ahler manifold. The key remark is 
that 
because of \eqref{doi}, $\Phi^{-1}(0)$ is the cone $N'\times \RR_+$ over 
$N'=\{x\in N\; ;\; (x,1)\in \Phi^{-1}(0)\}$. Moreover, since the actions of 
$G$ and $\rho_t$ commute, one has an induced action of $G$ on $N'$. Then  
$$\Phi^{-1}(0)/G\cong (N'\times \RR_+)/G\cong N'/G\times \RR_+$$
The manifold $N'/\!/G\times \RR_+$ is K\"ahler, as reduction of a K\"ahler manifold, 
but we still have to check that 
this K\"ahler structure is a cone one. 
For the more general, symplectic case, 
this was done in \cite{Bo-Ga 4}.
Let $g_0$ be the reduced K\"ahler metric 
and $g'$ be the Sasakian reduced metric on $N'/\!/G$.  
It is easily seen that the lift of $g_0$ to $\Phi^{-1}(0)$ coincides with the 
lift of the cone 
metric $r^2g'+dr^2$ on horizontal fields. 
This implies that the cone metric coincides with $g_0$.

Summing up we have proved:
\begin{te}
Let $(N,g,\xi)$ be a Sasakian manifold and let $(C(N),$ $C(g)$, $J)$ be the K\"ahler cone over it. 
Let a compact Lie group $G$ act by holomorphic isometries on $C(N)$ and 
commuting with the action of the $1$-parameter group generated by the field 
$R_0$. If a moment map with regular value $0$ exists for this action, then a 
moment map with regular value $0$ exists also for the induced action of $G$ on 
$N$. Moreover, the reduced space  $C(N)/\!/G$ is the K\"ahler cone  over the 
reduced Sasakian manifold $N/\!/G$.
\end{te}
The advantage of defining the Sasakian reduction \emph{via} K\"ahler reduction, 
as done in \cite{Bo-Ga-Ma 1} for $3$-Sasakian manifolds, is the avoiding of 
 curvature computations.

\section{Examples: $S^1$ actions on Sasakian spheres}
\begin{ex}
Start with $S^7\subset \CC^4$ with its standard Sasakian structure. Let the 
complex coordinates of $\CC^4$ be $(z_0,...,z_3)$, with $z_j=x_j+iy_j$.
The 
contact form on $S^7$ can then be written 
$$\eta=\sum_{j=0}^3(x_jdy_j-y_jdx_j)$$
and its Reeb field is 
$$\xi=\sum_{j=0}^3(x_j\partial y_j-y_j\partial x_j).$$
Let $S^1$ act on $S^7$ by $e^{it}\mapsto(e^{-it}z_0,e^{-it}z_1,e^{it}z_2,e^{it}z_3)$. 
The associated field of this action is (in real coordinates) 
\begin{equation*}
\begin{split}
X_0&=-(x_0\partial y_0-y_0\partial x_0)-(x_1\partial y_1-y_1\partial x_1)+\\
&+(x_2\partial y_2-y_2\partial x_2)+(x_3\partial y_3-y_3\partial x_3).
\end{split}
\end{equation*}
The moment map $\mu:S^7\rightarrow \RR$ reads:
$$\mu(z)=\eta_z(X_0)=-\abs{z_0}^2-\abs{z_1}^2+\abs{z_2}^2+\abs{z_3}^2$$
with zero level set 
$$\{z\in S^7\; ;\; \abs{z_0}^2+\abs{z_1}^2=\abs{z_2}^2+\abs{z_3}^2\}=
S^3(\frac{1}{\sqrt{2}})\times S^3(\frac{1}{\sqrt{2}}).$$
Clearly $\mu$ is nondegenerate on $\mu^{-1}(0)$.

    The reduced space can be identified with $S^3\times S^3/S^1$ which, by 
\cite{W-Z}, is diffeomorphic with $S^2\times S^3$. (In this case, one can also avoid 
the topological arguments in \cite{W-Z} and identify the reduced 
space by observing that the following diffeomorphism of $S^3\times S^3$: 
$(z_0,z_1,z_2,z_3)\mapsto (z_1z_4+\overline{z_2z_3}, z_1z_3-\overline{z_2z_4}, z_3,z_4)$ 
is equivariant with respect to the previous $S^1$ action which restricted to the second 
factor of the product is the usual action inducing the Hopf fibration; \emph{mille grazie} 
to Rosa Gini and  Maurizio Parton for letting us know it, \cite{GP}).

    The reduced Sasakian structure obtained in this way on $S^2\times S^3$ is 
easily checked to be Einstein and to project on the K\"ahler Einstein metric of 
$\CC P^1\times \CC P^1$ making the fibre map be a Riemannian submersion. 
As by \cite{W-Z} such an Einstein metric is unique, our reduced Sasakian structure 
coincides with the Sasakian structure found in \cite{Or-Pi 3} viewing $S^2\times S^3$ 
as minimal submanifold of $S^7$, total space of the pull-back over  
$\CC P^1\times \CC P^1$ of the Hopf bundle $S^7\rightarrow \CC P^3$. The same  
Einstein-Sasakian metric on $S^2\times S^3$ also appears in \cite{Ta}, constructed 
by a different approach.
\end{ex}

\begin{ex}
Consider again $S^7$ as starting Sasakian manifold, but let $S^1$ act by:
$e^{it}\mapsto(e^{-kit}z_0,e^{it}z_1,e^{it}z_2,e^{it}z_3)$, $k\in\ZZ_+$.
Now $\mu^{-1}(0)\cong S^1(\sqrt{\frac{k}{k+1}})\times
S^5(\sqrt{\frac{1}{k+1}})$. 
In order to identify the reduced space, consider the $k:1$ mapping
$$S^1\times S^5\ni (z_0,z_1,z_2,z_3)\mapsto ((z_0)^{-k},z_1,z_2,z_3)\in 
S^1\times S^5.$$
It induces a $k:1$ map from $M=S^1\times S^5/S^1$, where 
$S^1$ acts 
diagonally, to 
the reduced space $\mu^{-1}(0)/S^1$ with the action given above. 
As in \cite{GP}, the map 
$$(z_0,...,z_3)\mapsto (z_0,
\overline{z_0}z_1,\overline{z_0}z_2,\overline{z_0}z_3)$$
is an equivariant diffeomorphism of $S^1\times S^5$, equivariant with
respect to the diagonal action of $S^1$ and the action of $s^1$ on the
first factor. Hence $M$ is diffeomorphic to $S^5$ and the reduced
Sasakian space is $S^5/\ZZ_k$.

\end{ex}
\begin{ex}
In general, consider the weighted action of $S^1$ on 
$S^{2n-1}\subset \CC^{n}$ by:
$$(e^{it},(z_0,...,z_{n-1}))\mapsto (e^{\la_0it}z_0,...,e^{\la_nit}z_{n-1})$$
where $(\la_0,...,\la_{n-1})\in \ZZ^{n}$. The associated moment map 
$$\mu(z)=\la_0\abs{z_0}^2+...+\la_n\abs{z_{n-1}}^2$$
is regular on $\mu^{-1}(0)$ for any  $(\la_0,...,\la_{n-1})$ such that 
$\la_0...\la_{n-1}\neq 0$, $(\la_0,...,\la_{n-1})=1$ and  at least 
two $\lambda$'s have different signs (compare with the $3$-Sasakian case where 
the weights obey to more restrictions, cf. \cite{Bo-Ga-Ma 1}).

    Now take $\la_0=...=\la_k=a$ and $\la_{k+1}=...=\la_{n-1}=-b$, $a,b\in \ZZ_+$ 
relatively prime. 
Then  $\mu^{-1}(0)\cong S^{2k+1}(\sqrt{\frac{a}{a+b}})\times 
S^{2(n-k)-1}(\sqrt{\frac{b}{a+b}})$. Note that the induced metric on $\mu^{-1}(0)$ 
coincides with the product metric of the standard metrics of the two
factors. We then see that the reduced space is diffeomorphic with an $S^1$
factor of the above product of spheres given by the following action:
$$(e^{it},(x,y))\mapsto (e^{iat}x,e^{-ibt}y).$$
One can now adapt the arguments of \cite{W-Z}, Cor. 2.2 and prove that the
reduced spaces are $S^1$ bundles over $\CC P^k\times \CC P^{n-k-1}$ and,
for $1\leq k$, $4<n$, they are not homeomorphic to each other in general.

	However, for $k=1$, $n=4$, the reduced space is always diffeomorphic with
$S^2\times S^3$. Hence, one obtains an infinite family of 
Sasakian structures  on $S^2\times S^3$.

	Note also that if $n$ is \emph{even},  choosing like in the first example, the
first half of the $\la$'s to be $-1$, the rest of them $1$, the reduced Sasakian metric
is \emph{Einstein}, again acoording to \cite{W-Z}. 
 \end{ex}

\end{document}